\documentclass{article}
\usepackage[paper=a4paper,dvips,top=3cm,left=3cm,right=3cm,
   foot=1.4cm,bottom=3cm]{geometry}
\usepackage{amsfonts, amsmath}
\usepackage{amsthm}
\usepackage{color}
\usepackage{abstract}
\usepackage{fixltx2e}
\usepackage{algorithm}
\usepackage{algpseudocode}
\usepackage{graphicx}
\usepackage{cite}
\usepackage{listings}

\newcommand{\vt}[1]{{\bf #1}}
\newcommand{\x}{\vt{x}}

\newcommand{\Ten}{\mathcal{T}}
\newcommand{\REAL}{\mathbb{R}}

\newcommand{\R}{\mathbb{R}}

\newcommand{\T}{\top}
\newcommand{\ts}{\mathcal}
\newtheorem{lemma}{Lemma}[section]

\newcommand\keywords[1]{\textbf{Keywords}: #1}
\newcommand\subclass[1]{\textbf{AMS subject classifications}: #1}
\usepackage{fancyhdr}

\newtheorem{theorem}{Theorem}[section]
\newtheorem{definition}{Definition}[section]
\newtheorem{Proposition}{Proposition}[section]

\begin{document}
\title{
A gradient projection method for  semi-supervised hypergraph clustering problems}

\author{
  {Jingya Chang}\thanks{%
    School of Mathematics and Statistics, Guangdong University of Technology
    ({\tt jychang@zzu.edu.cn}). This author's work was partially supported by the National Natural Science Foundation of China (No. 11901118 and 62073087).}
 \and
  {Dongdong Liu}\thanks{%
    School of Mathematics and Statistics, Guangdong University of Technology
    ({\tt ddliu@gdut.edu.cn}). This author's work was partially supported by National Natural Science Foundation of China (No. 12101136), Guangdong Basic and Applied Basic Research Foundations (Nos.  2023A1515011633, 2020A1515110967), Project of Science and Technology of Guangzhou (No. 202102020273).}
 \and
  {Min Xi}\thanks{Corresponding author. School of Mathematics and Statistics, Guangdong University of Foreign Studies
    ({\tt mxi@gdufs.edu.cn}).  This author's work was partially supported by the National Natural Science Foundation
of China (No. 12171168), Guangdong Province Higher Education Foundation (No. 2021ZDZX1071).}
}
\date{}

\maketitle

\begin{abstract}
Semi-supervised clustering problems focus on clustering data with labels. In this paper, we consider the  semi-supervised hypergraph problems. We use the hypergraph related tensor to construct  an orthogonal constrained optimization model. The optimization problem is solved by a retraction method, which employs the polar decomposition to map the gradient direction in the tangent space to the Stefiel manifold. A nonmonotone curvilinear  search  is implemented to guarantee reduction in the objective function value. Convergence analysis demonstrates that the first order optimality condition is satisfied at the accumulation point. Experiments on synthetic hypergraph and hypergraph given by real data demonstrate the effectivity of our method.
\end{abstract}


\subclass {05C65, 15A69, 65K05, 90C35}
\section{Introduction}
Clustering and classification are two important tasks in machine learning. Clustering approaches aim to divide a number of items without labels into several groups,  while classification methods provide a classifier with the help of  labeled data and classify other data by using the classifier. On one hand data labeling  in real life, such as  getting  labels in computer-aided diagnosis or part-of-speech tagging,  is usually time-consuming or difficult \cite{van2020survey}. On the other hand, sometimes in reality, few annotated  points are melded in the unannotated data set in  clustering problems and  taking advantage of the priori label information  often enhances the clustering performances.

Semi-supervised learning is to complete the learning task based on both the labeled and unlabeled data \cite{zhou2018brief}. Semi-supervised clustering approach has wide applications in different areas. In image processing, the semi-supervised clustering approach was employed for image classification and segmentation \cite{enguehard2019semi}. A semi-supervised algorithm was proposed in \cite{kushagra2019semi} to solve the data de-duplication problem. For microarray expression data analysis, the knowledge from gene ontology data set was well utilized to generate the clustering algorithm \cite{cheng2004knowledge,fang2006knowledge}. An active semi-supervised clustering method was applied to modeling the complex industrial process \cite{lei2016modeling}.

Hypergraph is a useful tool to save and describe the high dimensional and complex  data arising from reality \cite{zhou2022clustering,veldt2022hypergraph}. Therefore, we consider the semi-clustering problems that are modeled by hypergraphs. We employ the multi way array to generate the clustering costs. By relaxing the label value and imposing the label matrix on the Stefiel manifold, we construct a tensor related optimization model and utilize the projection operator as a retraction to compute the feasible descent direction on the manifold. Numerical experiments show that our method works well on synthetic and real data.

The outline of this paper is as follows. In Section 2, we introduce the preliminary knowledge. The semi-clustering optimization model is given in Section 3, while the computing algorithm and convergent result are presented in Section 4 and Section 5 respectively. The numerical performance of our method is demonstrated in Section 6. Finally, we conclude our work in Section 7.

\section{Preliminaries}
In this section we demonstrate some useful notions and   results on hypergraphs and tensors.
Let $\R^{[r,n]}$ be the $r$th order $n$-dimensional real-valued tensor space, i.e.,
$$\R^{[r,n]}\equiv \R^{\overbrace{{n \times n \times \cdots \times n}}^{r\text{-times}}}.$$
The  tensor $\Ten=(t_{i_1 \cdots i_r}) \in\REAL^{[r,n]}$  with $ i_j=1,\ldots,n$ for  $j=1,\ldots,r,$ is said to be symmetric if $t_{i_1 \cdots i_r}$ is unchanged under any permutation of  indices \cite{chen2016positive}. Two operations  between
$\Ten$ and any vector $\x \in \mathbb{R}^n$ are stipulated as
\begin{equation*}
    \Ten\x^r \equiv \sum_{i_1=1}^n\cdots\sum_{i_r=1}^n
      t_{i_1\cdots i_r}\x_{i_1}\cdots \x_{i_r}
\end{equation*}
and \begin{equation*}
    (\Ten\x^{r-1})_i \equiv \sum_{i_2=1}^n\cdots\sum_{i_r=1}^n
      t_{ii_2\cdots i_r}\x_{i_2}\cdots \x_{i_r}, \quad \, \text{for} \quad i=1,\ldots,n.
\end{equation*}
Note that,  $\Ten\x^r \in \mathbb{R}$ and $\Ten\x^{r-1} \in \mathbb{R}^n$  are a scalar and a vector respectively,  and $\Ten\x^r = \x^
{\T} (\Ten\x^{r-1}).$ The tensor outer product of $\mathcal{A}\in\mathbb{R}^{[p,n]}$ and $\mathcal{B}\in\mathbb{R}^{[q,n]}$ is given by
\[\mathcal{A}\circ\mathcal{B}=\left(a_{i_1i_2\ldots i_p}b_{j_1j_2\ldots j_q}\right)\in\mathbb{R}^{[p+q,n]}. \]

\begin{definition}[Hypergraph]
A  hypergraph is defined as $G=(V,E)$, where $V=\{1,2,\ldots,n\}$ is the vertex set and
  $E=\{e_1,e_2,\ldots,e_m\}\subseteq 2^{V}$ (the powerset of $V$) is the edge set. We call $G$ an $r$-uniform hypergraph when  $|e_p|=r \geq 2$ for $p=1,\ldots,m$ and $e_i \neq e_j$ in case of $i \neq j.$

  If each edge of a hypergraph is linked with a positive number $s(e),$ then this hyperpragh is called a weighted hypergraph and $s(e)$ is the weight associated with the edge $e.$   An ordinary hypergraph  can be regarded as a weighted hypergraph with the weight of  each edge being $1.$
\end{definition}

The Stiefel manifold is $\mathcal{M}_n^p:=\{X: X^TX=I\}.$ We take the Euclidean metric as the Riemann metric on the Stiefel manifold and its tangent space.

\section{Semi-supervised Clustering Model}
Consider the hypergraph semi-supervised clustering problem. Our task is to cluster the $n$ vertices of the hypergraph into $k$ groups according to the hypergraph structure, while few categorization labels are given. Here $k$ is the number of clusters  and is usually much less than $n.$

Denote an indicator matrix $X\in R^{n\times k}$ as
\begin{equation}\label{IndicatorMatrix}
  X_{ij}=\left\{
           \begin{array}{ll}
             1, & \hbox{the $i$th vertex is in the $j$th cluster;} \\
             0, & \hbox{otherwise.}
           \end{array}
         \right.
\end{equation}
 $X$ is a  matrix with $0$-$1$ elements and its columns are orthogonal to each other. Thus, the clustering cost of an edge $e$ whose vertices are not labeled \cite{chang2020hypergraph} can be described as
\begin{equation*}
\sum_{j=1}^k  s_e \sum_{ i_1,i_2\in e} \left|\frac{X_{i_1j}}{\sqrt[r]{d_{i_1}}}-\frac{X_{i_2j}}{\sqrt[r]{d_{i_2}}}\right|^r.
\end{equation*}
The symbols $d_{i_1}$ and $d_{i_2}$ are degrees of vertices $i_1$ and $i_2$ respectively, and $r$ is the order of the hypergraph.
The total cutting cost of the  hypergraph  is
\begin{equation}\label{Cost_Function_1}
 f(X) =  \sum_{j=1}^k \sum_{e\in E} s_e \sum_{ i_1,i_2\in e} \left|\frac{X_{i_1j}}{\sqrt[r]{d_{i_1}}}-\frac{X_{i_2j}}{\sqrt[r]{d_{i_2}}}\right|^r.
\end{equation}
It is shown in \cite[Proposition 3.1]{chang2020hypergraph} that the cost function is equivalent to
  \begin{equation}\label{Cost_Function_1a}
    f(X) = \sum_{j=1}^{k}\ts{L} \vt{x}_j^r.
  \end{equation}
Here $$\mathcal{L} =  \sum_{e\in E} s_e \sum_{i,j\in e}
     \underbrace{\vt{u}_{ij}\circ \vt{u}_{ij}\circ \cdots \circ \vt{u}_{ij}}_{r \ \ \text{times} } $$ is a tensor related with the weighted hypergraph $G = (V,E,s),$
 $\vt{u}_{ij}=\frac{\vt{e}_i}{\sqrt[r]{d_i}}-\frac{\vt{e}_j}{\sqrt[r]{d_j}}$ with $\vt{e}_i$ and $\vt{e}_j$ being the $i$th  and $j$th columns of the identity matrix respectively, and $\vt{x}_j = X(:,j).$

%

Next, we take into account the vertices labeled. Denote a matrix $Y$  to save the labels of items as follows
\begin{equation}\label{MatrixG}
  Y_{ij}=\left\{
         \begin{array}{ll}
           1, & \hbox{the $i$th item is known in the $j$th cluster;} \\
           0, & \hbox{otherwise.}
         \end{array}
       \right.
\end{equation}
We call the matrix $Y$ the label matrix herein. It is natural that we try  best to retain the labeled vertex in its predetermined cluster during the clustering process. Therefore,  we add a regularization term $\lambda\|X_Y-Y\|_F^2$   to keep our clustering results insistence with the given labels. Here
$$ \quad
 (X_Y)_{ij}=\left\{
         \begin{array}{ll}
           X_{ij}, & \hbox{if $Y_{ij}\neq0$} \\
           0, & \hbox{otherwise.}
         \end{array}
       \right.$$
  Hence, the cutting cost becomes $\sum_i \mathcal{L}x_i^r+\lambda\|X_Y-Y\|_F^2$ under the semi-supervised condition. In terms of the $0$-$1$ constraint of the elements in $X,$ we relax it to the orthogonal constraint $ X^T X=I.$ The question is then transformed from discrete to continuous, and is easier to handle than the original $0$-$1$ one. Finally we get our orthogonal constrained model of the semi-supervised clustering problem as follows
\begin{equation}\label{Model}
\left\{
\begin{aligned}
&\min f(X)=\sum_{i=1}^k \mathcal{L}\vt{x}_i^r+\lambda\|X_Y-Y\|_F^2  \\
&s.t. \  X^T X=I.
\end{aligned}
\right.
\end{equation}

\section{Computation}
In this section, we first review  the gradient of a function on the tangent space, and then we introduce the algorithm based on this gradient direction.
%

For optimization models constrained on the Steifel manifold, there are two threads to follow \cite{hu2020brief, jiang2015framework, absil2009optimization}. One way is transforming the constrained optimization problem into an unconstrained one by
using mathematical programming techniques such as the penalty method, the augmented Lagrangian method \cite{xiao2021solving} and then solve it by unconstrained optimization methods. The other way is first finding a descent direction in the tangent space of the current point, then mapping an appropriate point in the descent direction to the Stiefel manifold.

We adopt the  second route to  compute the orthogonal constrained problem \eqref{Model}. The tangent space at a point $X\in S_{n,p}$ is $$T_X := \{Z:X^T Z+Z^T X=0\}.$$  Suppose function $f(X):\mathcal{M}_n^p \rightarrow R$ is differentiable. A retraction is a smooth mapping from the tangent bundle to the manifold \cite{absil2009optimization,wen2013feasible,
chang2020hypergraph}. Then for any $ X \in \mathcal{M}_n^p,$ the retraction $h_X(Z)$ is a mapping from $T_X\in R^{n\times p}$ to $ \mathcal{M}_n^p$ with $h_X(0)=X.$ The function $f(h_X(Z)):T_X \rightarrow R$ is also differentiable. Denote the gradient of function $f$ at the current point $X$ as $G.$  Because the objective function $f(X)$ is separable, we can compute the gradient $G$ in parallel. Take the inner product of two matrices $\langle A,B \rangle$ as $\langle A,B \rangle=tr(A^T B).$ For any vector $\xi$ in the tangent space of $T_X(M)$ the  gradient of $f$ at $X$ projected onto the tangent space satisfies $$\langle \nabla f(X),\xi\rangle=Df(x)[\xi].$$ In this paper, we take  the direction $\nabla f(X)$  as $ G-X\frac{X^TG+G^T X}{2}.$ Next we show that  $\nabla f(X)=0$ for $X^TX=I$  is equivalent to the  first order optimality condition of the constrained model \eqref{Model}.
\begin{Proposition}\label{lemmaOptimalityCondition}
Assume $X\in S_{n,p}.$   Let $\nabla f(X)= G-X\frac{X^TG+G^T X}{2}.$  The first order optimality conditions of \eqref{Model} hold if and only if
\begin{equation}\label{nablaf=0}
\nabla f(X)=0,
\end{equation}
with $X^{\top} X=I.$
\end{Proposition}

\begin{proof}
Let the symmetric matrix $\Lambda$ be the Lagrangian multiplier of the constraint $X^{\top}X=I.$  The Lagrangian function is
$f(X)+\frac{1}{2}tr(\Lambda^{\top}(X^{\top}X-I)).$ Similar to the first order optimization conditions in \cite{wen2013feasible}, we get the first order optimality of \eqref{Model}
\begin{equation}\label{kkt}
         G-XG^{\top}X=0, \quad \quad X^{\top}X=I
      \end{equation}
with   $ \Lambda = X^{\top}G.$

If the first order optimality conditions \eqref{kkt} hold, we have $X^{\top}G=XG^{\top}.$ Then
$$\nabla f=G-X\frac{X^TG+G^T X}{2}=G-XG^{\top}X=0.$$ On the other hand if $\nabla f =0,$  we can also get $X^{\top}G=XG^{\top},$ which validates  the first optimality conditions \eqref{kkt}.
\end{proof}

The direction $-\nabla f$  is employed for the descent direction in the tangent space. Our second step is  mapping the decent direction $-\nabla f$ to the manifold by using a retraction. We  project $-\nabla f$   to the manifold by utilizing the polar decomposition \cite{absil2009optimization,manton2002optimization}. The projection of any vector $Z$ onto the Stefiel manifold is defined as $$h(Z)= \underset{Q\in \mathcal{M}_n^p }{\arg\min} \|Z-Q\|_F^2.$$ If the SVD of $Z$ is $Z=U\Sigma V^T$, then the optimal solution of $h(Z)$ can be computed by $UV^T$ \cite{golub2013matrix}. Also, $h(Z)$ can be expressed as $Z(Z^TZ)^{-\frac{1}{2}} $ equivalently.

\begin{Proposition}\label{lemma4.2}
Suppose $X\in \mathcal{M}_n^p$ and  $Z$ is a vector in the tangent space $T_X.$ Consider the  univariate function $$h_X^Z(t)=\underset{Q\in \mathcal{M}_n^p }{\arg\min} \|Q-(X+tZ)\|_F^2$$ with its domain $t\in R,$ which is the projection of $X+tZ$ onto the Stiefel manifold.
 The derivative of $h_X^Z(t)$ at $t=0$ is
\begin{equation}\label{h'0}
(h_X^Z)'(0)=Z.
\end{equation}
When $Z=-\nabla f(X),$ the derivative of $f(h_X^Z(t))$ at $t=0$ is
\begin{equation}
f'_t(h_X^Z(0))=-\|\nabla f(X)\|^2.
\end{equation}
\end{Proposition}

\begin{proof}
For any arbitrary matrix $\tilde{Z}\in R^{n\times p},$ it can be decomposed as
$\tilde{Z}=XA+X_\bot B+XC,$ in which $A\in R^{p\times p}$ is skew-symmetric and $C\in R^{p\times p}$ is symmetric. It is proved in \cite[Lemma 8]{manton2002optimization} that,
\begin{equation}\label{MantonL8}
h_X^{\tilde{Z}}(t)=X+t(XA+X_\bot B)+O(t^2).
\end{equation}
On the other hand, the tangent vector space at $X$ is $T_X(M)=\{X\Omega+X_\bot K, \Omega=-\Omega, K\in R^{(n-p)\times p}\}$  \cite{absil2009optimization}. Therefore, for $Z\in T_X(M)$
\begin{equation}\label{lemma1}
h_X^{Z}(t)=X+tZ+O(t^2).
\end{equation}
Since $h_X^{Z}(0)=X,$ we get $(h_X^Z)'(0)=Z.$

By the chain rule, when $Z=-\nabla f(X)$ the derivative of $f(h_X^Z(t))$ with respect to $t$ at $t=0$ is
\begin{eqnarray*}\label{chainrule}
  f_t^{'}(h_X^Z(0)) &=& \langle f'(X),(h_X^Z)'(0)\rangle=\langle G,-\nabla f(X)\rangle \\
    &=&  -\langle \nabla f(X),\nabla f(X)\rangle-\frac{1}{2}\langle \nabla f(X),X(X^T G+G^T X)\rangle \\
&=& -\langle \nabla f(X),\nabla f(X)\rangle-\frac{1}{2}[\langle X^T\nabla f(X),X^TG+G^TX\rangle] \\
&=& -\langle \nabla f(X),\nabla f(X)\rangle-\frac{1}{4}[\langle X^TG-G^TX,X^TG+G^TX\rangle]\\
&=& -\|\nabla f(X)\|^2.
\end{eqnarray*}
\end{proof}

We employ a nonmonotone line search method  to find a proper step size along $X+tZ.$  Given $X_k$ and the descent direction $Z=-\nabla f(X_k),$
by using the adaptive feasible BB-like method proposed in \cite{jiang2015framework}, we find a step size $t_k$ such that
\begin{equation}\label{Armijo0}
f(h_X^{Z}(t_k))\leq f_r+\delta t_k f'(h_X^{Z}(0)),
\end{equation}
where $f_r$ is a reference objective function value.
Let $S_{k-1}=X_k-X_{k-1}, Y_{k-1}=-\nabla f(X_k)+\nabla f(X_{k-1}).$ The parameter $t_k$ takes the following BB stepsizes alternately:
$$t_k^1=\frac{\langle S_{k-1}, S_{k-1}\rangle}{|\langle S_{k-1},Y_{k-1} \rangle|},t_k^2=\frac{|\langle S_{k-1},Y_{k-1} \rangle|}{\langle Y_{k-1}, Y_{k-1}\rangle}.$$ Let $L$ be a preassigned positive integer and $f_{best}$ be the  current best function value. Denote by $f_C$ the  maximum objective function value after  $f_{best}$ is found. The reference function value $f_r$ is updated only when the best function value is not improved in $L$ iterations. The detailed steps are shown below.


\begin{equation}\label{fr}
\begin{split}
& \texttt{if} \qquad f_{k+1}<f_{best}\\
& \hspace{3em} f_{best}=f_{k+1},f_c=f_{k+1},l=0\\
& \texttt{else}   \\
& \hspace{3em} f_c=\max\{f_{k+1},f_c\},l=l+1,\\
&   \hspace{3em} \text{if} \ l=L, f_r=f_c,f_c=f_{k+1},l=0 \ \text{end} \\
& \texttt{end}
\end{split}
\end{equation}

 The proposed method is shown in Algorithm \ref{Algorithm}.

\begin{algorithm}
\caption{SSHC: Semi-supervised Hypergraph Clustering Algorithm}
\begin{algorithmic}[1]
\Require $0<\rho<1, X_0\in \mathcal{M}_n^p,0<\beta<1,\epsilon>0,L$ be a positive integer.
\While{}
\State Compute $G$ at the point $X_k$ in parallel.
\State Compute $\nabla f(X_k), Z = -\nabla f(X_k).$
\State Find the smallest nonnegative integer $m$, such that
\begin{equation}\label{ArmijoAlg}
f(h_{X_k}^{Z}(\beta^m \hat{t_k}))\leq f_r+\delta\beta^m \hat{t_k} f'(h_{X_k}^{Z}(0)).
\end{equation}
holds and set $t_k=\beta^m \hat{t_k}$ with $\hat{t_k}=t_k^1$ or $t_k^2$  alternately.
\State Update $X_{k+1}$ and $f_r, f_{best},f_c $ by \eqref{fr}.
\EndWhile
\end{algorithmic}\label{Algorithm}
\end{algorithm}

\section{Convergence Analysis}
The sequence $\{\nabla f(X_k)\} $ generated from Algorithm \ref{Algorithm} either terminates with $\nabla f(X_k)=0$ or is infinite.   By reductio and absurdum,  we prove that when the iteration is infinite, a subsequence of $\{\nabla f(X_k)\} $ converges to 0,
 which means
$$\underset{k\rightarrow \infty}{\lim \textmd{inf}}\|\nabla f(X_k)\|=0.$$
\begin{lemma}\label{lemma2}
Assume  there exists a constant $\varepsilon>0$ such that
\begin{equation}\label{assumption}
\|\nabla f(X_k)\|\geq \varepsilon.
\end{equation}
Under this assumption, the step size $t_k$ generated by \eqref{ArmijoAlg}  satisfies
\begin{equation}\label{stepsizeBound}
t_k\geq c,
\end{equation}
where $c$ is a constant.
\end{lemma}
\begin{proof}
Suppose the conclusion does not hold. Then we can find a subsequence $\{k_i\}$ satisfies that
$$t_{k_i}\rightarrow 0\ \ \text{as} \ k_i \rightarrow \infty.$$
We use the symbol $k$ instead of $k_i$ for simplicity. According to Proposition \ref{lemma4.2}, when $Z=-\nabla f(X_k)$ the Taylor expansion of $f(h^Z_{X_k}(t))$ at the point $t=0$ is
\begin{eqnarray}\label{taylorExp}
  f(h^Z_{X_k}(t)) &=& f(h^Z_{X_k}(0))+t f'_t(h^Z_{X_k}(0))+o(t) \nonumber \\
   &=& f(X_k)-t\|\nabla f(X_k)\|^2  +o(t).
\end{eqnarray}
If $\hat{t_k}$ is not accepted in the Armijo-type search \eqref{ArmijoAlg}, then we have $\hat{t_k}\geq \beta^{-1}t_k$ and
\begin{eqnarray}\label{ineq1Lemma}
   f(h^Z_{X_k}(\beta^{-1}t_k)) &>& f_r+\delta \beta^{-1}t_k f'(h^Z_{X_k}(0) \nonumber \\
   &\geq& f(X_k)-\delta \beta^{-1}t_k \|\nabla f(X_k)\|^2.
\end{eqnarray}
Substituting $\beta^{-1}t_k$ for $t$ in \eqref{taylorExp}, we have
\begin{equation}\label{ineq2Lemma}
  f(h^Z_{X_k}(\beta^{-1}t_k)) = f(X_k)-\beta^{-1}t_k\|\nabla f(X_k)\|^2  +o(t_k).
\end{equation}
Since $t_{k}\rightarrow 0\  \text{as} \ k \rightarrow \infty,$ by combining \eqref{ineq1Lemma} and \eqref{ineq2Lemma}, we obtain
$$(1-\delta)\|\nabla f(X_k)\|^2+o(1)\leq 0.$$
This inequity is impossible, when $\|\nabla f(X_k)\|\geq \varepsilon.$ The proof is then completed.
\end{proof}

\begin{theorem}
If the sequence $\{X_k\}$ given by Algorithm \ref{Algorithm} is infinite, then we have
\begin{equation}\label{theorem}
\underset{k\rightarrow \infty}{\lim \inf}\|\nabla f(X_k)\|=0.
\end{equation}
\end{theorem}
\begin{proof}
Suppose  this conclusion is not true, then the assumption \eqref{assumption} holds. We save all values of  $f_r$ in \eqref{fr} in the sequence
$\{f_r^{m}\},$ where the index $m$ means the $m$th value of $f_r.$ Denote the index of the first iteration that is produced from the line search \eqref{Armijo0} related with $f_r^{m}$ as $k_m.$  Let $l_m$ be the index number that satisfies $f(X_{l_m})=\max\limits_{k_m\leq j< k_{m+1}}{f(X_j)}.$  From  \eqref{ArmijoAlg}, we have
\begin{equation}\label{ThEq1}
 f(X_{l_m})\leq f_r^{m}-\delta t_{l_m} \|\nabla f(X_{l_m})\|^2.
\end{equation}
Also from the updating process we have
\begin{equation}\label{ThEq2}
f_r^{m+1}\leq f(X_{l_m}).
\end{equation}
By \eqref{ThEq1}, \eqref{ThEq2} and \eqref{stepsizeBound}, we obtain
\begin{equation}\label{ThEq3}
f_r^{m+1}\leq f_r^{m}-\delta c \|\nabla  f(X_{l_m})\|^2.
\end{equation}
If $\{f_r^m\}$ is finite, the sequence of $\{f_{best}\}$ is infinite which contradicts with the fact that  $f(X)$ is bounded below. Therefore,  $\{f_r^m\}$ is an infinite sequence. Then based on Lemma \ref{lemma2} and \eqref{ThEq3} we get
\begin{equation}\label{ThEq4}
+\infty>\sum_{m=1}^{+\infty}[f_r^{m}-f_r^{m+1}] \geq \sum_{m=1}^{+\infty} \delta c \|\nabla f(X_{l_m})\|^2,
\end{equation}
which indicates that the assumption \eqref{assumption} is impossible. The conclusion \eqref{theorem} is finally proved.
\end{proof}

\section{Numerical Experiments}
In this section, we demonstrate the numerical performance of SSHC method for clustering synthetic and real data. For each problem, we run 100 times and report the maximum, minimum, average and median accuracy rate by recording the average value of the results given by the 100 runs.

\subsection{Semi-supervised and supervised clustering of artificial hypergraphs}
\begin{figure}
  \centering
  \includegraphics[width=0.35\textwidth]{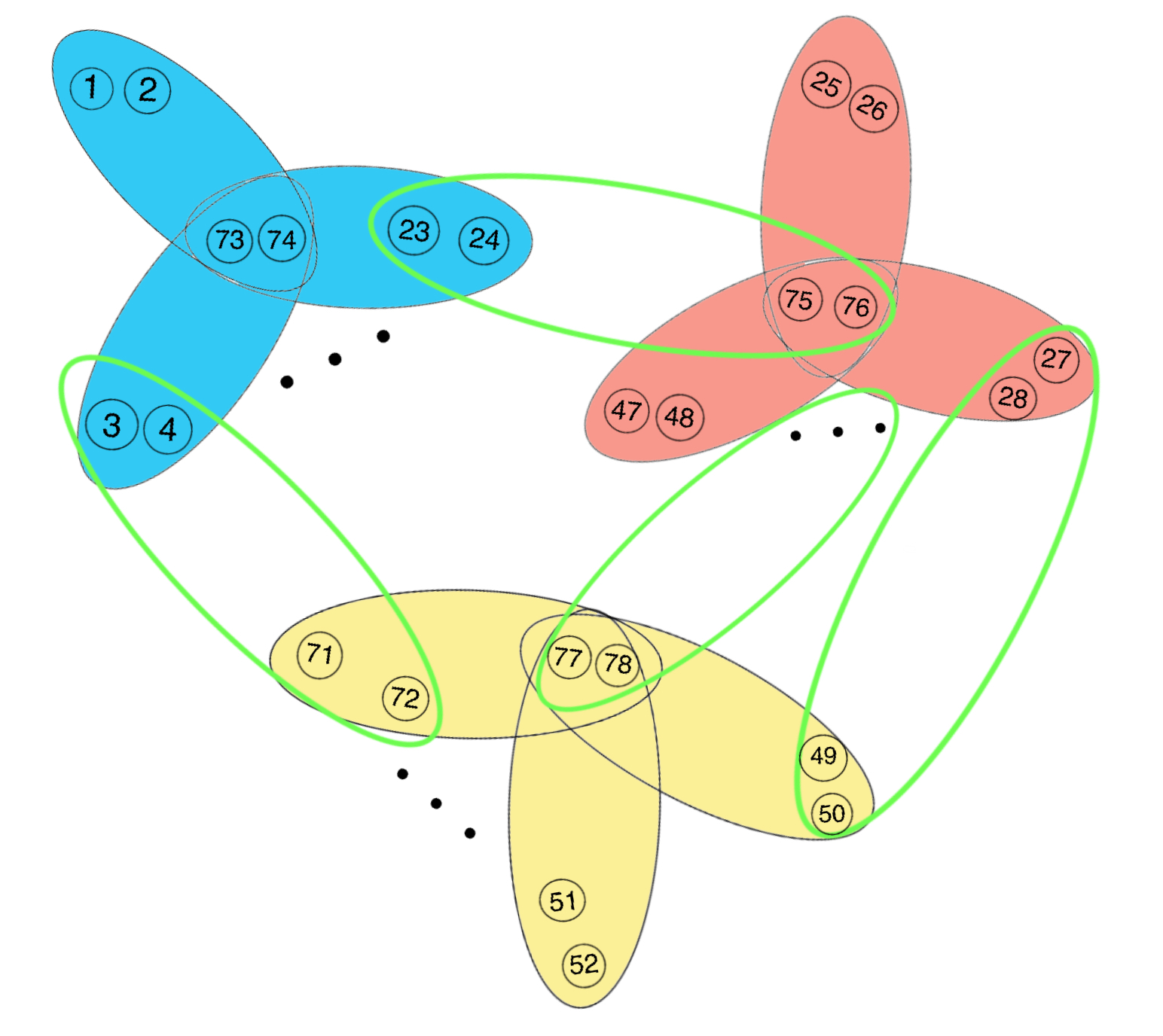}
  \caption{The $4$-uniform hyerpgraph }\label{4graph}
\end{figure}
In this subsection, we employ the proposed SSHC method to cluster a class of  artificial hypergraphs. We compare the SSHC method and the method SHC, which is in fact the unsupervised model by replacing the objective function in \eqref{Model} as
$$f(X)=\sum_{i=1}^k \mathcal{L}\vt{x}_i^r.$$    In order to construct a $4$-uniform hyerpgraph, we first generate 3 $4$-uniform sub-hypergraphs. Each of the sub-hypergraph has $12$ edges which share two common vertices. Then we randomly produce 4 more edges that contain vertices from vertex sets of different sub-hypergraphs. An example of the final hypergraph is shown in  Figure \ref{4graph}.  The weight of this hypergraph is an all one vector. The clustering results of this hypergraph is demonstrated in Table \ref{SyntheticData}. For the semi-supervised problem , 10 percent of the vertices are labeled. It can be seen that the clustering accuracy is promoted by the proposed method.

\begin{table}\caption{The clustering results with different proportions of labels}\label{SyntheticData}

\begin{center}
\begin{tabular}{c|cccc}
  \hline
    & max& min & average & median \\ \hline
  SSHC & 0.145833 &       0.062500   &    0.095000    &   0.093750     \\
  SHC & 0.125000   &     0.125000  &     0.125000   &    0.125000    \\
  \hline
\end{tabular}
\end{center}
\end{table}

\subsection{Yale face data clustering}
\begin{figure}
  \centering
  \includegraphics[width=0.7\textwidth]{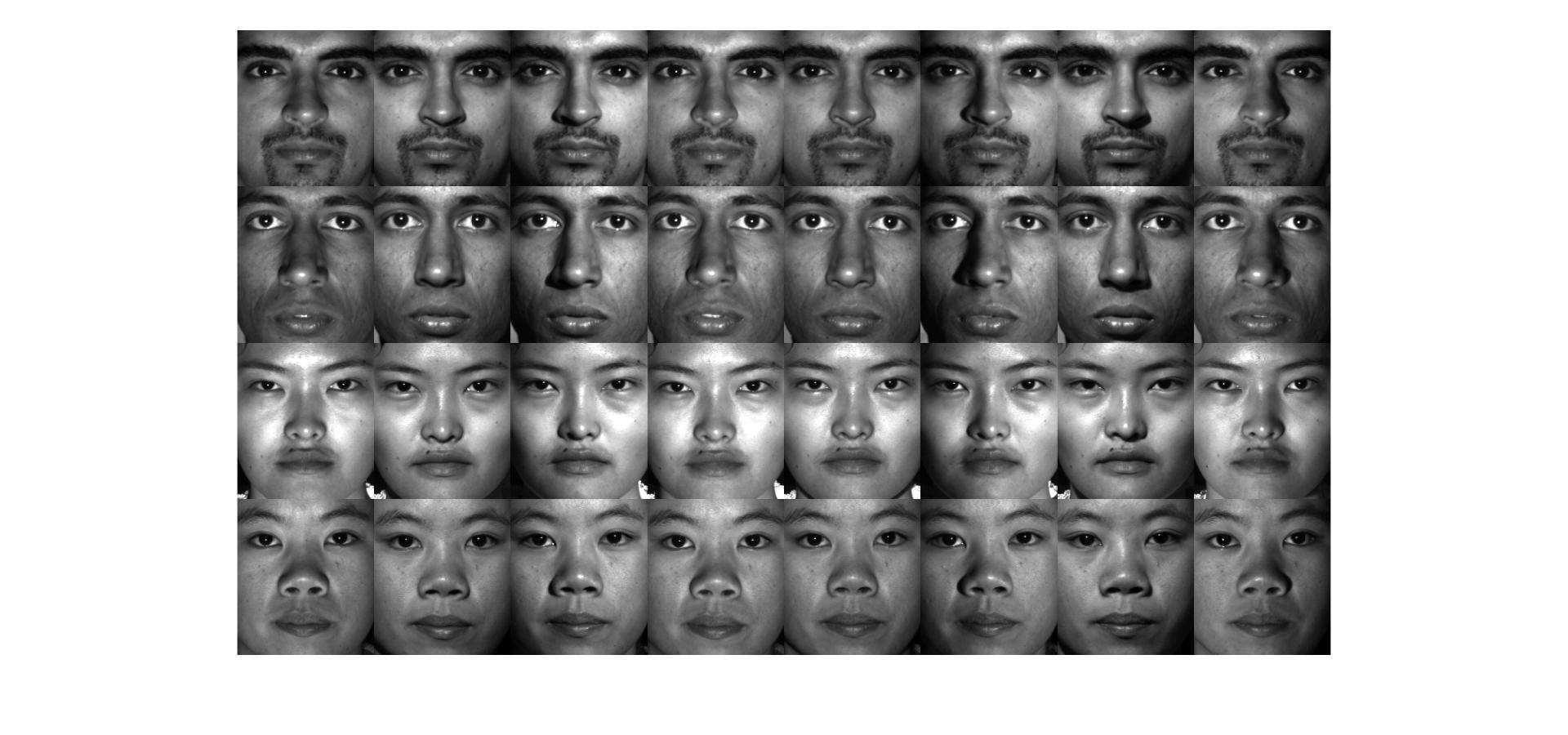}
  \caption{8 images of 4 person from the Yale data set}\label{FigYaleface}
\end{figure}
The extended Yale face data base  B contains 600 face images of 30 persons under 20 lighting conditions \cite{georghiades2001few,lee2005acquiring}. Before computation, each image is resized into $ 48\times 42$ pixels and expressed as a vector. Figure \ref{FigYaleface} displays 8 images of 4 persons as examples. Our task is to group the images of each person into a cluster. The hypergraph is constructed based on a rough clustering result by k-means method. The images that are in one cluster are put in one edge. Also, we randomly choose other images combing with the images in a cluster to produce more edges. The weights of  edges are given based on the Euclidian distance of  its image vectors.

We cluster the images for 100 times and record the average error, median error, max error and minimum error. The ratio of labeled images to all is shown in the $Ratio$ column.  From Table \ref{Yaleface} we can see that, the clustering error decreases when the number of labeled images increase.
\begin{table}\caption{The clustering results with different proportions of labels}\label{Yaleface}
\begin{tabular}{c|cccc|c|cccc}
  \hline
  Ratio &average & max& min & median &Ratio & average &max& min & median  \\ \hline
       0 & 0.036   &   0.375  &  0   & 0   & 0.1   &  0.022  &  0.343    &   0&0   \\
  0.2    &  0.019   &   0.375  &  0& 0     &0.3    &  0.013   &  0.344   & 0&0   \\
   \hline
\end{tabular}
\end{table}
\section{Conclusion}
In this paper, we give a tensor related optimization model to compute the hypergraph clustering problems with little part of labels provided.
We use the  polar decomposition as a retraction on the Stiefel manifold. The convergence analysis shows that an accumulation point of the iteration sequence is a stationary point. Numerical experiments indicate that the method improves the computation accuracy when compared to the unsupervised model. However, the effectiveness of the hypregraph clustering method relies on an appropriate hypergraph of the data. The construction of a hypergraph that reasonably reflects the data structure and relationship is a meaningful topic in our future research.

 \bibliographystyle{abbrv}
\bibliography{bibfile}

\end{document}